\documentclass[11pt]{amsart}

\usepackage{graphicx}
\usepackage{amsmath}%
\usepackage{amsfonts}%
\usepackage{amssymb}

\usepackage{epsfig}

\usepackage{psfrag}

\newtheorem{theorem}{Theorem}[section]

\newtheorem{definition}[theorem]{Definition}

\newtheorem{lemma}[theorem]{Lemma}

\newtheorem{proposition}[theorem]{Proposition}

\numberwithin{equation}{section}

\def\N{{\mathbb{N}}}

\def\lam{\lambda}
\def\Th{\theta}
\def\gam{\gamma}
\def\Om{\Omega}
\def\om{\omega}

\def\const{{\rm const}}

\def\nula{\nu_\lam}
\def\wt{\widetilde}
\def\dint{\int\!\!\!\int}
\def\Diag{{\rm Diag}}

    \newcommand{\Leb}{{\mathcal L}} 
\def\Zk{{\mathcal Z}}
\def\Bk{{\mathcal B}}

\newcommand{\Ek}{{\mathcal E}}

\newcommand{\eps}{{\varepsilon}} 
 
\newcommand{\half}{{\frac{1}{2}}}
\def\Ik{{\mathcal I}}

\begin{document}

\title[Spacings for finite Bernoulli convolutions]{Spacings and pair correlations for finite Bernoulli convolutions}

\author{Itai Benjamini}
\address{Itai Benjamini,
Department of Mathematics, Weizmann Institute of Science,
Rehovot, 76100, Israel}
\author{Boris Solomyak}
\address{Boris Solomyak, Box 354350, Department of Mathematics,
University of Washington, Seattle WA 98195}
\email{solomyak@math.washington.edu}

\thanks{ 
\indent Research of Solomyak was
partially supported by NSF grant \#DMS-0654408.} 

\begin{abstract}

We consider finite Bernoulli convolutions with a parameter $1/2 < \lam < 1$ supported on a
discrete point set, generically  of size $2^N$. These sequences are uniformly distributed with respect to the infinite Bernoulli convolution
measure $\nula$, as $N\to \infty$.
Numerical evidence
suggests that for a
generic $\lam$, the distribution of spacings between appropriately rescaled points is Poissonian.
We obtain some partial results in this direction; for
instance, we show that, on average, the pair correlations do not exhibit attraction or repulsion in the limit. On the other hand, for certain
algebraic $\lam$ the behavior is totally different.

\end{abstract}


\maketitle

\thispagestyle{empty}

\section{Introduction}

With original motivation coming from physics \cite{BT}, many authors investigated spacings and other statistical properties of various number-theoretic
sequences,  see e.g.\ \cite{marklof,rusa,ruza,rusaza} 
and references therein. Given a sequence $\{\Th_n\}$ equidistributed in the unit interval, the  nearest-neighbor spacing distribution is obtained
by ordering the first elements of the sequence $\Th_{1,N} \le \Th_{2,N} \le \cdots \le \Th_{N,N}$, and then defining normalized spacings to be
$$
\delta_n^{(N)}:= N (\Th_{n+1,N} - \Th_{n,N}).
$$
One can also consider next-to-nearest neighbor spacings or more generally, for any fixed $\ell\ge 1$,
$$
\delta_{\ell,n}^{(N)}:= N (\Th_{n+\ell,N} - \Th_{n,N}).
$$
A function $P_\ell(s)$ is called the limiting distribution function of $\{\delta_{\ell,n}^{(N)}\}$ as $N\to \infty$ if, for any interval $[a,b] \subset 
[0,1)$,
\begin{equation} \label{eq-level.ell}
\lim_{N\to\infty} \frac{1}{N} \#\bigl\{n\le N:\ \delta_{n,\ell}^{(N)} \in [a,b]\bigr\} = \int_a^b P_\ell(s)\,ds.
\end{equation}
If the points $\Th_n$ are independent uniformly distributed on $[0,1]$ random variables, we get the Poisson model, with
$$
P_\ell(s) = \frac{s^{(\ell-1)}}{(\ell-1)!} e^{-s},\ \ \mbox{for all}\ \ell\ge 1.
$$
Other widely used statistical measures are correlation functions.  The pair correlation
function for the sequence $\{\Th_n\}$ is defined by
$$
R_2(s, \{\Th_n\}, N) = \frac{1}{N}  \#\Bigl\{(n,m):\ n\ne m, \ n,m\le N,\ |\Th_{n,N} - \Th_{m,N}| \le s \Bigr\}.
$$
In the Poisson model we have $R_2(s, \{\Th_n\}, N) \to 2s$, as $N\to \infty$.

In this paper, we initiate a study of a similar problem for sequences arising from overlapping iterated function systems and non-standard digit
expansions.
For $\lam \in (\half,1)$ let
$$
A_N(\lam) = \Bigl\{(1-\lam)\sum_{n=0}^{N-1} a_n \lam^n:\ a_n\in \{0,1\}\Bigr\}.
$$
We are interested in the randomness properties of these sets as $N\to \infty$. Note that a point in 
$A_N(\lam)$ may have more than one representation; in other words, we may have ``coincidences.'' This will happen
whenever $\lam$ is a zero of a polynomial with coefficients $\{-1,0,1\}$, which we call $\{0,\pm 1\}$ polynomials for short. 
In this case we take the point ``with multiplicity,'' so that our
sequence  always has size $2^N$. Consider the uniform measure $\nu_\lam^N$ on $A_N(\lam)$; this is a discrete
measure, with a point in $A_N(\lam)$ having mass $2^{-N}$ times the number of representations. It is well-known that for all
$\lam$,
the measures $\nu_\lam^N$ converge vaguely to the infinite Bernoulli convolution measure $\nula$, defined as the 
distribution of the random series $\sum_{n=0}^\infty \om_n (1-\lam)\lam^n$, where $\om_n$ are i.i.d.\ Bernoulli random
variables taking the values $0,1$ with equal probability. Often $\pm 1$ are used as digits instead of $0,1$; this is just
a linear change of variable. We also multiplied by $(1-\lam)$, so that the measures are all supported on $[0,1]$. Note 
that $(\half,1)$ is the interesting parameter range for this problem: for $\lam< 1/2$ the measure $\nula$ is supported on a
Cantor set, and $A_N(1/2)$ is just the set of binary rationals in $[0,1)$ with denominators less than $2^N$; also, for
$\lam \le 1/2$ the ordering of $A_N(\lam)$ agrees with the lexicographical ordering on $\{0,1\}^N$.

The measures $\nula$, for $\lam \in (\half,1)$, 
have been much studied, but are still rather mysterious, with many unsolved problems,
see \cite{sixty} and references therein. 
Besides mathematics, they come up in physics, 
information theory, and economics, see e.g.\ \cite{KR,sidor,econ}.
Briefly, it is known that $\nula$ is continuous for all $\lam$, it is either purely absolutely continuous or purely singular,
it is absolutely continuous for a.e.\ $\lam \in  (\half,1)$, and it is singular for Pisot reciprocals (see the definition
below). It is an open problem whether these are the only exceptions.
The only explicitly known $\lam$ for which $\nula$ is absolutely continuous are Garsia reciprocals (see the definition
below), which include $2^{-1/k}$, for $k \ge 2$. Moreover, the roots of $1/2$ are the only parameters for which the density is explicitly computable.
(When we talk about ``density of $\nula$'' in this paper, we mean the Radon-Nikodym derivative $\frac{d\nula}{dx}$.)

Let $\{\xi_{j,N}\}_{j=1}^{2^N}$ be the ordering of $A_N(\lam)$, counting with multiplicity, so that
$$0 = \xi_{1,N} \le \xi_{2,N} \le \cdots \le \xi_{2^N,N}= 1-\lam^N.
$$
Spacing statistics provide a measure of randomness for the sequences $\{\xi_{n,N}\}$ approximating $\nula$.
Numerical evidence seems to indicate that, for a typical $\lam$, there exist level spacing distributions. 
However, we first need to rescale the sequence. Let 
$$
\wt{\xi}_{n,N} = F_\lam(\xi_{n,N}),\ \ n=1,\ldots,2^N,\ \ \ \mbox{where}\ F_\lam(\xi) = \nula(-\infty,\xi).
$$
These sequences are uniformly distributed modulo one, in the sense that for any interval $[a,b]
\subset [0,1]$, the number of points $\wt{\xi}_{n,N}$ in $[a,b]$ divided by $2^N$ tends to $b-a$. 
(This follows from the fact that $\nula$ is continuous, hence the cumulative distribution function (CDF)
of $\nula^N$ converges to the CDF of $\nula$ everywhere.)

\medskip

{\bf Conjecture 1.} {\em For almost every $\lam\in (\half,1)$ the rescaled sequences $\wt{\xi}_{n,N}$ 
are distributed according to the Poisson model.}

\medskip

We have numerical evidence in support of this conjecture,
however, proving it seems far beyond our reach. A key difficulty is that very few measures $\nula$ are known explicitly. 
Therefore, we mostly focus on the pair correlation function for the non-rescaled sequences
$$
R_2(s, \lam, 2^N) :=
\frac{1}{2^N}  \#\Bigl\{(x,y):\ x\ne y, \ x,y\in A_N(\lam),
\ |x-y| \le \frac{s}{2^N} \Bigr\}
$$
which provide information about the variance of average spacings.
A more modest goal than proving Conjecture 1
would be to verify the following

\medskip

{\bf Conjecture 2.} {\em For almost every $\lam\in (\half,1)$, there exist $c, C>0$ such that 
$$
cs \le R_2(s,\lam, 2^N) \le Cs\ \ \mbox{for all $N$ and $s>0$}.
$$
}

\medskip

This would mean that our sequences exhibit {\bf no attraction} and {\bf no repulsion}.
Our main results (Theorems 2.1 and 2.3) makes a step in this direction, but we have to do averaging over the parameter $\lam$.

Note that for certain $\lam$ the behavior of the spacings is definitely not Poissonian. For example, let 
$\lam = 1/\Th$, where $\Th$ is a Pisot number, that is, an algebraic integer $>1$ whose conjugates are  all less than one in
modulus. Such $\nula$, with $\lam\in (\half,1)$, are the only singular infinite Bernoulli convolutions known \cite{erdos1}. 
In this case, Garsia's Lemma \cite{garsia1} says that the spacings between distinct points in $A_N(\lam)$
are all bounded below by $C\cdot \lam^N$, which implies that there are ``massive'' coincidences, since 
$2^N - C\cdot \lam^N$ spacings are zeros. This is an extreme case of attraction.
Something similar happens for all zeros of $\{0,\pm 1\}$ polynomials.

On the other hand, consider the
case of ``Garsia numbers.'' We say that $\Th$ is Garsia if it is a zero of a monic polynomial with integer coefficients and
constant term $\pm 2$, whose conjugates are all greater than one in modulus.
Garsia reciprocals include $2^{-1/k}$ for $k\ge 2$, 
but also many other numbers (e.g.\ the reciprocal of the positive root of $x^3 -2x-2 =0$, which is about $.5652$, 
or the reciprocal of the positive root of $x^3-x^2-2=0$, which is about $.5898$).
Garsia reciprocals $\lam$
are the only {\em specific} (as opposed to generic) $\lam \in (\half,1)$ for which it is known that $\nula$ is 
absolutely continuous, and moreover, they have bounded density. Again by Garsia's Lemma \cite{garsia1}, it is known that the
spacings between points in $A_N(\lam)$ are bounded below by $C\cdot 2^{-N}$, and there are no coincidences. 
This is an extreme case of repulsion. 

It is interesting to compare
our situation with that for sequences $\alpha n^d$ mod 1, for irrational
$\alpha$ and $d\ge 2$.
It is conjectured that for $\alpha$ badly approximable by rationals the spacing distribution is Poissonian, and partial results in this direction are obtained
in \cite{rusa,rusaza}, in particular, for a.e.\ $\alpha$, the pair correlations converge to $2s$ as $N\to \infty$. On the other hand, \cite{rusa,rusaza} show that for
some irrational $\alpha$ the behavior of spacings is different.
In our problem, we prove that for $\lam$ which are extremely well-approximable by zeros of $\{0,\pm 1\}$ polynomials, the
pair correlation function $R_2(s,\lam,2^N)$ exhibits some attraction (Theorem 2.2).

There are several difficulties in pushing our results further; one of them
is that the distribution of zeros of $\{0,\pm 1\}$ polynomials is not well-understood.
In contrast with \cite{rusaza}, we cannot say anything about higher-order correlations; this is related to hard open problems on Bernoulli
convolutions. A technical assumption which we have to make in some results is that the parameter $\lam$ belongs to the ``transversality interval,''
that is, an interval where $\{0,\pm 1\}$ power series cannot have double zeros. 
Transversality \cite{posi,solo1,peso1,peso2} has been crucially important in the work on $\nula$.

We note that there are some features of our sequences which are ``non-random.'' In particular, some spacings tend to have very high multiplicity,
in other words, the distribution of non-rescaled spacings sometimes exhibit ``spikes.'' This is explained by the fact that
every spacing $g$ is a zero of a $\{0,\pm 1\}$ polynomial.
For every such  polynomial with $k$ zero coefficients, the same difference between points in $A_N(\lam)$ will be achieved $2^k$ times.
If $g$ is a very small number, it is likely to be a spacing, and if $N-k$ is small, we will see a ``spike."
This irregularity should dissappear after rescaling by the CDF of $\nula$.

Further, it turns out that the smallest spacing between points 
of $A_N(\lam)$, for a.e.\ $\lam$ in a transversality interval, is bounded below by $3^{-N} N^{-1-\eps}$, for infinitely many $N$, for any $\eps>0$ 
(follows from Theorem 2.4),
which is much larger than $(2^{-N})^2$ which we would get if the points were
sampled randomly and independently from the uniform distribution.

In Figures 1-3 we show the histogram of the distribution of nearest-neighbor spacings and also spacings with $\ell=2,3$
for the rescaled points $\{\wt{\xi}_{n,N}\}$, where 
$\lam=0.70880447$ was obtained from a random number generator, uniformly in the interval $[0.69,0.71]$. It was chosen so as to be reasonably close to 
$2^{-1/2} = 0.7017\ldots$ for which the CDF of the  infinite Bernoulli convolution measure is known precisely, and this CDF 
was used for the rescaling.
The figures show reasonable agreement with the Poisson distribution (indicated by the curves). Note that the units on the horizontal axis are
$2^{-N}$.
More details about the computations, {\em Mathematica} program, and additional simulations are
collected in the Appendix at the end of the paper. 

\begin{figure}[ht]
\begin{center}
\epsfxsize=4.5in
\epsffile{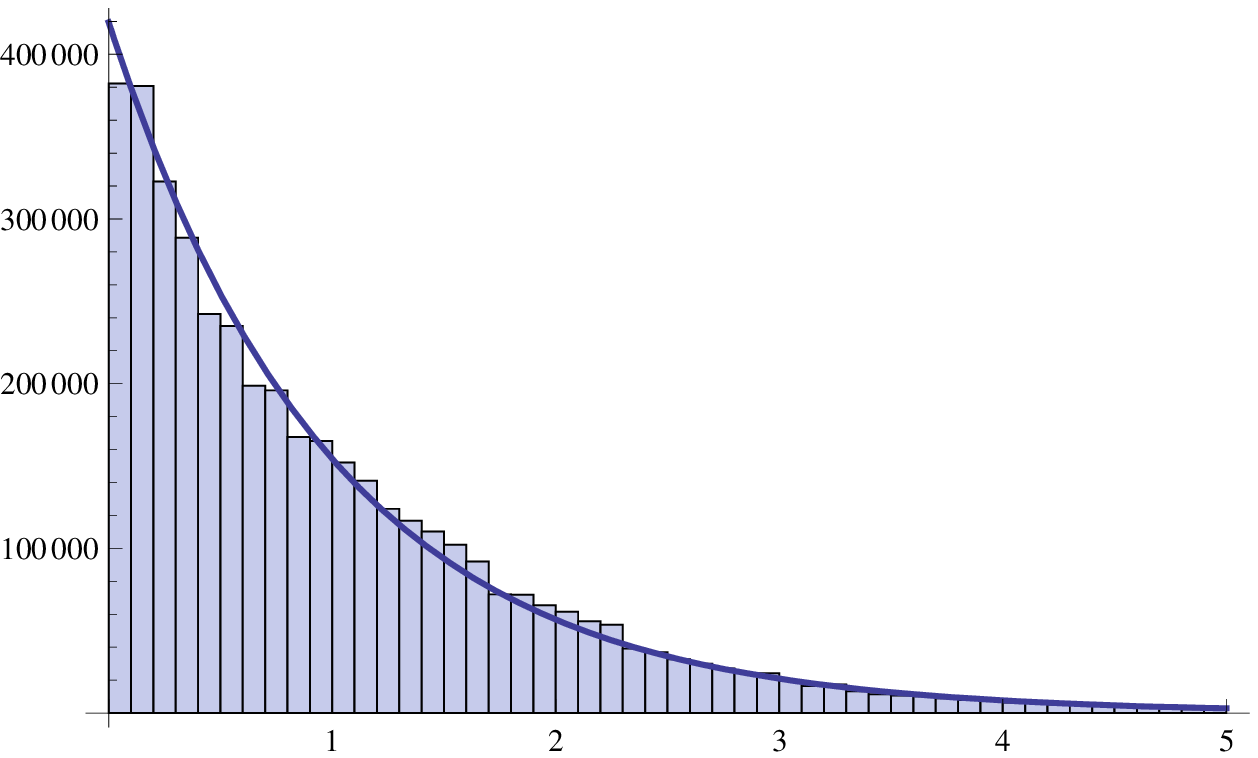}
\end{center}
\caption{Nearest-neighbor spacings for rescaled points, $\lam=.70880447,\ N=22$.}
\end{figure}

\begin{figure}[ht]
\begin{center}
\epsfxsize=4.5in
\epsffile{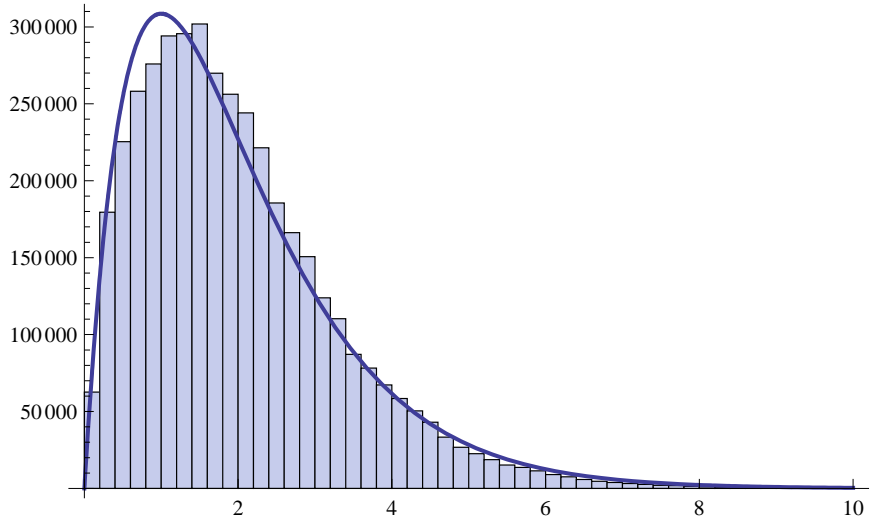}
\end{center}
\caption{Next-to-nearest-neighbor spacings for rescaled points, $\lam=.70880447,\ N=22$.}
\end{figure}

\begin{figure}[ht]
\begin{center}
\epsfxsize=4.5in
\epsffile{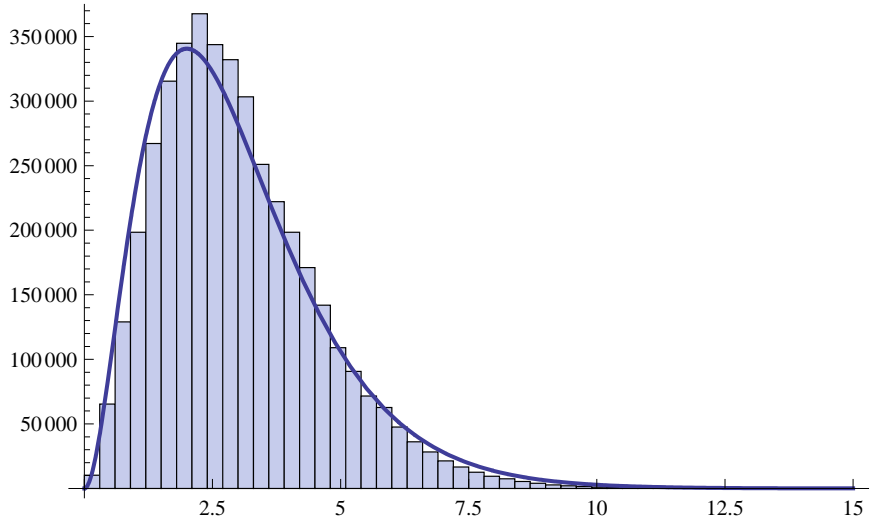}
\end{center} 
\caption{Spacings with $\ell=3$ for rescaled points, $\lam=.70880447,\ N=22$.}
\end{figure}

\section{Results} 

\begin{theorem} \label{th-noattr}
Let $I=[\lam_0,\lam_1]\subset (\half,0.668)$.
Then there exists $C_1$ (depending on $I$) such that
\begin{equation} \label{eq-noattr}
\forall\ N\ge 1,\ \forall\ s>0,\ \int_I R_2(s, \lam, 2^N) \,d\lam \le C_1s,
\end{equation}
\end{theorem}

{\bf Remarks.} 1.
The appearance of the number $0.668$ is due to ``transversality;'' it is explained in the next section. We do not know if the result holds all
the way to 1, although numerical evidence indicates that things actually get better as $\lam$ increases.

2. By Fatou's Lemma, we get $\int_I \liminf_{N\to \infty} R_2(s, \lam, 2^N) \,d\lam \le C_1s$ for $s>0$.
This shows that, on average, our spacing distributions exhibit no attraction, at least in the liminf sense.
We note that in other problems, for instance, concerning the distribution of $\{n\alpha\}$ for irrational $\alpha$, it is also
common to average over the parameter, see \cite[Th.3.2]{marklof}.

3. An application of Borel-Cantelli shows, for instance, that for any $\eps>0$, for a.e.\ $\lam\in (\half,0.668)$,
\begin{equation} \label{eq-asym1}
\forall\,s>0,\ R_2(s,\lam,2^N) \le s N^{1+\eps}\ \ \mbox{for $N$ sufficiently large}.
\end{equation}

4. An estimate analogous to (\ref{eq-noattr}) follows for the sets $A_N(\lam)\cap J$, for an arbitrary fixed interval $J\subset 
[0,1]$. Let $\#(J,N) = \#(A_N(\lam)\cap J)$ and set
\begin{equation} \label{eq-cor1}
R_2(s, \lam,J, N) :=
\frac{1}{\#(J,N)}  \#\Bigl\{(x,y):\ x\ne y, \ x,y\in A_N(\lam)\cap J,
\ |x-y| \le \frac{s|J|}{\#(J,N)} \Bigr\}.
\end{equation}
Then, obviously,
$$
R_2(s, \lam,J, N)\le R_2\Bigl( \frac{s|J|}{\#(J,N)},\lam,2^N\Bigr) \frac{2^N}{\#(J,N)}.
$$
\begin{sloppypar}
Note that $\lim_{N\to \infty} \#(J,N)/2^N = \nula(J)\in (0,1)$, so we get an analog of (\ref{eq-noattr}) for $R_2(s, \lam,J, N)$, with 
a constant $C_{1,J}\approx C_1|J|/\nula(J)^2$ for large $N$.
\end{sloppypar}

5. If $\nula$ has a density bounded below by $\eta>0$ on an interval $J=[a,b]\subset [0,1]$, then 
$$
\wt{R}_2(s,\lam,\wt{J},N) \le R_2 (s/\eta, \lam, J, N),
$$
where $\wt{R}_2(s,\lam,\wt{J},N)$ is the correlation sum analogous to (\ref{eq-cor1}) for the rescaled points $\wt{\xi}_{n,N}$ in 
$\wt{J} = [F_\lam(a), F_\lam(b)]$ (recall that $F_\lam$ is the CDF of $\nula$).
It is known that the density $\nula$ is bounded away from zero on $[\eps,1-\eps]$ for $\eps>0$, for a.e.\ $\lam> 2^{-1/2}$, and numerical evidence
suggests that this is true for a typical $\lam > (\sqrt{5}-1)/2\approx .618$ as well.

\medskip

The following theorem shows that (\ref{eq-asym1}) fails for many points.

\begin{theorem} \label{th-attr}
{\bf (i)} 
For any $\lam_0$ which is a zero of a $\{0,\pm 1\}$ polynomial, there exists $\rho = \rho(\lam_0)< 2$ such that
\begin{equation} \label{eq-growth}
|A_N(\lam_0)|\le C_2\rho^N,
\end{equation}
where $|A_N(\lam_0)|$ is the size of the set \underline{\em without multiplicities}. Thus,
\begin{equation} \label{eq-asym3}
R_2(0,\lam_0,2^N) \ge  C_2^{-1}(2/\rho)^N-1.
\end{equation}

{\bf (ii)}
For any interval $I \subset (\half,1)$ and $\eps>0$,
there is an uncountable set $\Ek_I\subset I$, such that 
for $\lam\in \Ek_I$
\begin{equation} \label{eq-asym2}
\forall\,s>0,\ R_2(s,\lam,2^N) \ge 2^{N^{1-\eps}}\ \ \mbox{for infinitely many $N$}.
\end{equation}
\end{theorem}

{\bf Remark.} One can ask how large is the uncountable set $\Ek_I$ in terms of Hausdorff measure and dimension. Our proof certainly produces
a set of Hausdorff dimension zero, but it shouldn't be hard to find a gauge function for which the Hausdorff measure is positive. We leave this to the
interested reader.

\medskip

Next we turn to the estimates from below. 

\begin{theorem} \label{th-norepul}
For any $\lam_0\in (\half,1)$, $\eps\in (0,1-\lam_0)$, and $L>0$, there exists $C_3 = C_3(\lam_0,\eps)>0$ such that
\begin{equation} \label{eq-norepul}
\forall\ N\ge 1,\ \forall\ s\in (0,L),\ \ \ \int_{\lam_0-\eps}^{\lam_0+\eps} R_2(s, \lam, 2^N) \,d\lam \ge C_3 s.
\end{equation}
\end{theorem}

{\bf Remarks.} 
1. An analog of (\ref{eq-norepul}) holds for $R_2 (s, \lam, J, N)$, for any interval $J \subset [0,1]$, with a constant $C_{3,J}$.

2. If $\nula$ has a density bounded above by $M$, then the pair correlation for the rescaled points satisfies
$$
\wt{R}_2(s, \lam, 2^N) \ge R_2(s/M, \lam, 2^N).
$$
It is known \cite{solo1} that $\nula$ has bounded density for a.e.\ $\lam\in (2^{-1/2},1)$.

\medskip

It is often more convenient (especially for the proofs) to work with 
\begin{equation} \label{eq-prim}
A'_N(\lam) = \Bigl\{\sum_{n=0}^{N-1} a_n \lam^n:\ a_n\in \{0,1\}\Bigr\} = (1-\lam)^{-1}A_N(\lam)
\end{equation}
and the corresponding correlation expression
$$
R'_2(s,\lam,2^N) =
\frac{1}{2^N}  \#\Bigl\{(x,y):\ x\ne y, \ x,y\in A'_N(\lam),
\ |x-y| \le \frac{s}{2^N} \Bigr\}.
$$
Since $(1-\lam)^{-1}$ is uniformly bounded  from zero and from infinity in all our theorems, proving results for $R'_2(s,\lam,2^N)$ will
immediately imply the desired statements. 

\medskip

Next we consider the {\bf smallest} and the {\bf largest} gaps (spacings) for $A'_N(\lam)$, which we denote $g_N(\lam)$ and $G_N(\lam)$
respectively. It is obvious that $g_N(\lam) \le (1-\lam)^{-1}2^{-N}$. For Garsia reciprocals $\lam$ we have $g_N(\lam) \ge c\cdot 2^{-N}$.
It seems likely  $g_N(\lam) =o(2^{-N})$ for a typical $\lam$, but we do not have a proof of this. The difficulty, again, is that we do not know
much about the distribution of the zeros of $\{0,\pm 1\}$-polynomials.

We have a lower bound for $g_N(\lam)$ for a
typical $\lam$ in the ``transversality interval.''

\begin{theorem} \label{th-mingap}
Let $\alpha_N>0$ be such that $\sum_N 3^N\alpha_N < \infty$.
Then for Lebesgue a.e.\ $\lam \in I = [\half,0.668]$ we have
$$
g_N(\lam) > \alpha_N\ \ \mbox{for infinitely many $N$}.
$$
\end{theorem}


\medskip

Concerning the largest gap, we do not have new results; we just note that questions related to
$G_N(\lam)$ have been 
studied by Erd\H{o}s, Jo\'o and Komornik (though they were
couched in different terms). It is easy to see that $G_N(\lam) \le \lam^{N-1}$ \cite[Th.4(a)]{EJK}
(this follows by the ``greedy algorithm''), and hence $G_N(\lam) = \lam^{N-1}$, since the ordered
sequence $A'_N(\lam)$ starts with $0, \lam^{N-1}.$ One can ask, however, whether this largest
gap is realized away from the endpoints, as $N\to
\infty$. The answer is ``yes'' for
all $\lam \in (\half,\lam_g)$, where $\lam_g = (\sqrt{5}-1)/2\approx 0.618\ldots$ is the golden ratio \cite[Th.4(b)]{EJK}, in which case, for
odd $N$, 
$$1+\lam^2+\lam^4+\cdots + \lam^{N-3},\  1+\lam^2+\lam^4+\cdots + \lam^{N-3}+\lam^{N-1}$$ are two consecutive points in $A'_N(\lam)$.
Erd\H{o}s, Jo\'o and Komornik \cite[Th.5]{EJK2} also showed that $G_N(\lam)$ is achieved at most $K$ times, independent of $N$, 
and not at all in $(\eps,(1-\lam)^{-1}-\eps)$ for any $\eps>0$, for sufficiently large $N$, for all $\lam \in ( 2^{-1/2},1)$ which are not zeros
of $\{0,\pm 1\}$ polynomials.
It is not known what happens between $\lam_g$ and $2^{-1/2}$, though it is reasonable to expect similar
behavior (gaps of $o(\lam^N)$, except near the endpoints),
at least for typical $\lam$.

\section{Proofs}

We need a result about $\{0,\pm 1\}$ polynomials.
Actually, it is convenient to state it for power series; of course this
includes polynomials as a special case.
Denote
$$
\Bk = \Bigl\{1 + \sum_{n=1}^\infty c_n x^n,\ c_n \in \{-1,0,1\}\Bigr\}.
$$

\begin{definition}
An interval $I\subset (0,1)$ is called an interval of
{\bf transversality} for $\Bk$ if
for every $f\in \Bk$,
we have
\begin{equation} \label{eq-trans}
x\in I,\ f(x) =0 \ \Rightarrow\ f'(x) \ne 0.
\end{equation}
\end{definition}

\begin{proposition}\cite{pablosol} \label{prop-trans}
The interval $[\half,0.668]$ is
an interval of transversality for $\Bk$.
\end{proposition}

Transversality is used in the following lemma.
The Lebesgue measure on the line will be denoted by $\Leb$.

\begin{lemma} (\cite[Lem.\,4.2]{peso2}, see also \cite{posi}) \label{lem-trans}
Let $I$ be a closed interval of transversality for $\Bk$ and $g\in \Bk$. Then there exists $C>0$, depending only on $I$, such that
\begin{equation} \label{eq-trans1}
\forall\,\rho>0,\ \Leb \{\lam\in I:\, |g(\lam)| \le \rho \} \le C \rho.
\end{equation}
\end{lemma}

Recall that we will be working with $A'_N(\lam)$ and $R'_2(s,\lam,2^N)$, see (\ref{eq-prim}).

\medskip

{\em Proof of Theorem~\ref{th-noattr}.} Our proof mimics the argument in \cite{peso1}, but we have to 
work with finite sums instead of power series.
Let $\Om_N = \{0,1\}^N$ and $\pi_\lam^N(\om) = \sum_{n=0}^{N-1} \om_n \lam^n$ for $\om = \om_0\ldots\om_{N-1} \in \Om_N$.
Note that $A'_N(\lam) = \pi_\lam^N(\Om_N)$. We will assume that $\pi_\lam$ is 1-to-1 on $\Om_N$, in other words, $\lam$ is not
a zero of a $\{0,\pm 1\}$ polynomial, which only excludes a countable set from $I$.
Further, let $\mu_N$ be the uniform (counting) measure on $\Om_N$, with each point having equal weight $2^{-N}$.
We have, by definition,
$$
R'_2(s, \lam, N) = 2^N (\mu_N\times \mu_N) \{(\om,\tau) \in \Om_N^{2}:\, \om\ne \tau,
\ |\pi_\lam^N(\om) - \pi_\lam^N(\tau)| \le s2^{-N}\}.
$$
By Fubini's Theorem,
\begin{equation} \label{eq-Fub}
\int_I R'_2(s, \lam, N) \,d\lam = 2^N
\dint_{\Om^2_N\setminus \Diag} \!\!\!\!\Leb\{\lam\in I:\, |\pi_\lam^N(\om) - \pi_\lam^N(\tau)| \le s2^{-N}\}\,
d \mu_N(\om)\,d\mu_N(\tau),
\end{equation}
where $\Diag = \{(\om,\om):\ \om \in \Om\}$.
We have
\begin{equation} \label{eq-eq1}
\pi_\lam^N(\om) - \pi_\lam^N(\tau) = \sum_{n=0}^{N-1} (\om_n-\tau_n) \lam^n = \pm \lam^k g(\lam),
\end{equation}
where $k = |\om\wedge\tau| = \min\{n:\,\om_n \ne \tau_n\}\le N-1$, and $g$ is a $\{0, \pm 1\}$ polynomial with constant term $1$, hence
$g \in \Bk$. 
By (\ref{eq-eq1}), $|\pi_\lam^N(\om) - \pi_\lam^N(\tau)| \le s2^{-N}$ implies $|g(\lam)| \le \lam_0^{-k} s 2^{-N}=:\rho$ for
$\lam \in I = [\lam_0,\lam_1]$. Applying (\ref{eq-trans1}) yields
$$
\Leb\{\lam\in I:\, |\pi_\lam^N(\om) - \pi_\lam^N(\tau)| \le s2^{-N}\} \le  C \lam_0^{-k}  s2^{-N},
$$
and we obtain from (\ref{eq-Fub}):
\begin{eqnarray*}
\int_I R'_2(s, \lam, N)  \,d\lam & \le &
C s \dint_{\Om^2_n\setminus \Diag} \lam_0^{-|\om\wedge\tau|}\,d \mu_N(\om)\,d\mu_N(\tau) \\
& = & C s \sum_{k=0}^{N-1} \lam_0^{-k}(\mu_N\times \mu_N) \{(\om,\tau):\,|\om\wedge\tau|=k\} \\
& = & C s \sum_{k=0}^{N-1} \lam_0^{-k} 2^{-k-1} <   \frac{C\lam_0\cdot s}{2\lam_0-1}\,,
\end{eqnarray*}
proving (\ref{eq-noattr}). \qed

\medskip

{\em Proof of Theorem~\ref{th-attr}.} (i) Let $\lam_0$ be a zero of a $\{0,\pm 1\}$ polynomial, that is,
$$
1 + \sum_{n=1}^k c_n \lam^n = 0,\ \ \mbox{with}\ c_n \in \{0,\pm 1\}.
$$
We can write $c_n = u_n - v_n$ for some $u_n, v_n \in \{0,1\}$, $n=1,\ldots,k$. It follows that
$1+\sum_{n=1}^k u_n \lam_0^n = \sum_{n=1}^k v_n \lam_0^n$, which means that every occurrence (not just at the beginning) of 
the ``block'' $(1 u_1\ldots u_k)$ in a sequence $a_0\ldots a_{N-1}$ representing a point in $A'_N(\lam_0)$,
can be replaced by $(0u_1\ldots u_k)$, without changing the value of the sum, in other words, yielding the same point in $A'_N(\lam_0)$.
It follows that 
$
|A'_N(\lam_0)|
$
is bounded above by the number of ``words'' of length $N$ in the shift of finite type consisting of all $0,1$ sequences which do not contain the
``forbidden'' block $(1 u_1\ldots u_k)$, see \cite[Ch.2]{LM}. Further, it is known (see \cite{LM}) that the number of words of length $N$ in a shift 
of finite type is not greater than $\const\cdot\rho^N$, where $\rho$ is the topological entropy. The entropy of the full shift on two symbols is 2, and
every ``forbidden word'' results in a drop of entropy, see \cite[Th.4.4.7]{LM}. Thus, $\rho< 2$, and (\ref{eq-growth}) is proved.
Now (\ref{eq-asym3}) follows from the Cauchy-Schwartz inequality: the $2^N$ points of $A_N(\lam)$ (counted with multiplicity) are grouped into
at most $\const\cdot\rho^N$ distinct values, and the sum of the squares of multiplicities is equal to $2^N(R_2(0,\lam_0,2^N)+1)$ (the ``$+1$'' 
comes from the diagonal).

\medskip

(ii) Let $\Zk$ denote the set of zeros of $\{0,\pm 1\}$ polynomials. Observe that $\Zk$ is dense in $(\half,1)$. 
Indeed, given
$\lam \in (\half,1)$, for any $k$, we can find $c_1=1, c_2,\ldots, c_k \in \{0,1\}$ such that $1-\sum_{n=1}^k c_n \lam^n \in [0,\lam^k)$
(e.g.\ by the ``greedy'' algorithm).
Then $p(x) = 1-\sum_{n=1}^k c_n x^n$ has derivative $p'(x)\le -1$ and hence $p$ has a zero in $[\lam,\lam+\lam^k)$.


We construct points $\lam \in I$, extremely well approximable by points in $\Zk$, to ensure that (\ref{eq-asym2}) holds.
We do this inductively. First find $\lam_0 \in \Zk \cap I$ and find $N_0$ so large that $R_2(0, \lam_0, 2^{N_0}) \ge 2^{N_0^{1-\eps}}$, which is
possible by (\ref{eq-growth}).
Then we find a closed interval $I_1$ around $\lam_0$, so small that 
$$
R_2(1, \lam, 2^{N_0}) \ge 2^{N_0^{1-\eps}}\ \ \mbox{for}\ \lam\in I_1,
$$
which is possible by continuity. We continue this process and find closed intervals
$I\supset I_1 \supset \ldots \supset I_k \supset \ldots$ and positive integers $N_0 < N_1 < \ldots < N_k < \ldots $, so
that 
$$
R_2(2^{-k}, \lam, 2^{N_k}) \ge 2^{N_k^{1-\eps}}\ \ \mbox{for}\ \lam \in I_k,\ k\ge 1.
$$
We are using density of $\Zk$ and (\ref{eq-growth}) at each step (however, the exponential rate depends on the point in $\Zk$, so $N_{k+1}-N_k$
may be quite large).
Clearly $\lam \in \bigcap_{k=1}^\infty I_k$ will satisfy (\ref{eq-asym2}). Modifying the construction slightly, namely, choosing two
distinct points in $\Zk$ in an interval and two disjoint intervals around them at each stage, we obtain a Cantor set of points in $I$ for which
(\ref{eq-asym2}) holds. \qed

\medskip

{\em Proof of Theorem~\ref{th-norepul}.}
The claim (\ref{eq-norepul}) will follow if we prove that
\begin{equation} \label{eq-Ik}
\Ik(s) := \int_{\lam_0-\eps}^{\lam_0+\eps} R'_2(s, \lam, 2^N) \,d\lam \ge C'_3 s.
\end{equation}

For any $k\ge 1$ we can find $c_1=1, c_2, \ldots, c_k \in \{0,1\}$ by the greedy algorithm, so that
$$
1 - \sum_{n=1}^k c_n \lam_0^n \in [0, \lam_0^k).
$$
Choose $k$ so that $\lam_0^k < \eps/10$. We have $k\ge 4$, since $\lam_0>\half$.
The idea is that for some integer $\ell$, for $N$ sufficiently large, the pair of points $x(\lam), y(\lam) \in A'_N(\lam)$, where
$$
x(\lam) = 1 + \sum_{n=k+\ell+1}^{N-1} a_n \lam^n,\ \ \
y(\lam) = \sum_{n=1}^k c_n \lam^n + \sum_{n=k+\ell+1}^{N-1} b_n \lam^n,
$$
will satisfy $|x(\lam)- y(\lam)| < s/2^N$ for all $\lam$ in an interval of size $\ge \const\cdot s/2^N$ contained in $(\lam_0- \eps, \lam_0+\eps)$,
for all choices of $a_n, b_n \in \{0,1\}$. This will contribute $\ge \const\cdot (s/2^N) \cdot (4^{N-k-\ell-1}/2^N)= \const\cdot s$
to the integral in (\ref{eq-Ik}), proving the estimate. Now for the details.
Consider
\begin{equation} \label{eq-f}
f(\lam):= x(\lam)-y(\lam)= 1-c_1\lam - c_2\lam^2 - \ldots - c_k \lam^k + \sum_{n=k+\ell+1}^{N-1} (a_n - b_n) \lam^n.
\end{equation}
We have
$$
|f(\lam_0)| < \lam_0^k + \frac{\lam_0^{k+\ell+1}}{1-\lam_0} < \frac{\eps}{10}+\frac{\eps}{10}\,\frac{\lam_0^{\ell+1}}{(1-\lam_0)}\,.
$$
Assume that $\ell$ satisfies
$$
(1-\lam_0)^{-1} \lam_0^{\ell+1} <1
$$
(there will be another condition on $\ell$ below). Then 
\begin{equation} \label{eq-new1} 
|f(\lam_0)| < \frac{\eps}{5}\,.
\end{equation}
Now we need to estimate $f'(\lam)$ for $\lam \in (\lam_0-\eps,\lam_0+\eps)$. We have
\begin{equation} \label{eq-new2}
|f'(\lam)| \le \sum_{n=1}^\infty n\lam^{n-1} = (1-\lam)^{-2}.
\end{equation}
On the other hand, since $c_1=1$ and $c_n \ge 0$, we obtain from (\ref{eq-f}):
$$
f'(\lam) \le -1 + \sum_{n=k+\ell+1}^N n\lam^{n-1} \le -1 + \frac{k+\ell+2}{(1-\lam)^2} \lam^{k+\ell}.
$$
Let $\lam_1 = \lam_0 + \eps$ and assume that $\ell\in \N$ is such that 
$$\frac{k+\ell+2}{(1-\lam_1)^2} \,\lam_1^{k+\ell} \le \half\,,$$ then
\begin{equation} \label{eq-new3} 
f'(\lam) \le -1/2\ \ \mbox{for all}\ \lam\in (\lam_0-\eps,\lam_0+\eps).
\end{equation}
From (\ref{eq-new1}) and (\ref{eq-new3}) it follows that there exists $\gam \in [\lam_0-\frac{2\eps}{5}, \lam_0 + \frac{2\eps}{5}]$
such that $f(\gam) = 0$. Then (\ref{eq-new2}) implies that 
$$
|f(\lam)| < s2^{-N}\ \ \mbox{for}\ \lam\in (\gam - (1-\lam_1)^2 s 2^{-N}, \gam + (1-\lam_1)^2 s 2^{-N}),
$$
and we have $s2^{-N}(1-\lam_1)^2 < \frac{3\eps}{5}$ for large $N$. Thus, any choice of $a_n, b_n\in \{0,1\}$, with  $n\in
[k+\ell+1,N-1]\cap \N$, contributes at least $s2^{-N}(1-\lam_1)^2\cdot 2^{-N}$ to the integral $\Ik(s)$, and taking all
the $4^{N-k-\ell-1}$ choices into account yields (\ref{eq-norepul}). \qed

\medskip

{\em Proof of Theorem~\ref{th-mingap}.}
It is clear that $g_N(\lam) = \sum_{n=0}^{N-1} (a_n - b_n)\lam^n$
for some $a_n,b_n\in \{0,1\}$, hence
$g_N(\lam) = \lam^k |p_{N-k}(\lam)|$ for some polynomial $p_{N-k} \in \Bk$
of degree $N-k-1$, where
$k\in \{0,1,\ldots,N-2\}$.
If $g_N(\lam) \le \alpha_N$, then
$$0 \le |p_{N-k}(\lam)| \le  \alpha_N \lam^{-k} \le \alpha_N 2^k.$$
By Lemma~\ref{lem-trans},
$$
\Leb\{\lam\in I:\, |p_{N-k}(\lam)| \le \alpha_N 2^k \} \le C \alpha_N 2^k.
$$ 
There are $3^{N-k-1}$ polynomials of degree $N-k-1$ in $\Bk$. Let
$$
E_N = \Bigl\{\lam\in I:\ g_N(\lam) \le \alpha_N\Bigr\}.
$$
It follows that
$$
\Leb(E_N) \le 
C \alpha_N \sum_{k=0}^{N-1} 2^k\cdot 3^{N-k-1} < C\alpha_N 3^N.
$$
Recalling  that $\sum_N \alpha_N 3^N< \infty$ we conclude by Borel-Cantelli that
a.e.\ $\lam \in I$ does not
belong to $E_N$ for infinitely many $N$, proving the claim. \qed

\section{Appendix: numerical experiments}

In Figure 4 we show the histograms for the non-rescaled and rescaled data sets for $\lam = 0.6429$, which was chosen randomly.
For all rescalings we used the density for the parameter $2^{-1/2}$ given by
$$
F(x) = \left\{ \begin{array}{ll} ax^2/2, & x \in [0, b] \\
                                      a b (x-b), & x\in [b, 1-b] \\
                                      1 - a(1-x)^2/2, & x\in [1-b,1] \end{array} \right.,\ \ \mbox{where}
\ a = 1.5\sqrt{2} + 2,\ b = \sqrt{2}-1.
$$
This ``improved'' the distribution of spacings even for $\lam$ not very close to $2^{-1/2}$.
In the histogram, we used the
50 ``bins'' from $0.1\cdot 2^{-N}$ to $5\cdot 2^N$.

\begin{figure}[ht]
\begin{center}
$\begin{array}{cc}
\multicolumn{1}{l}{} &
        \multicolumn{1}{l}{} \\
\epsfxsize=3.5in
\epsffile{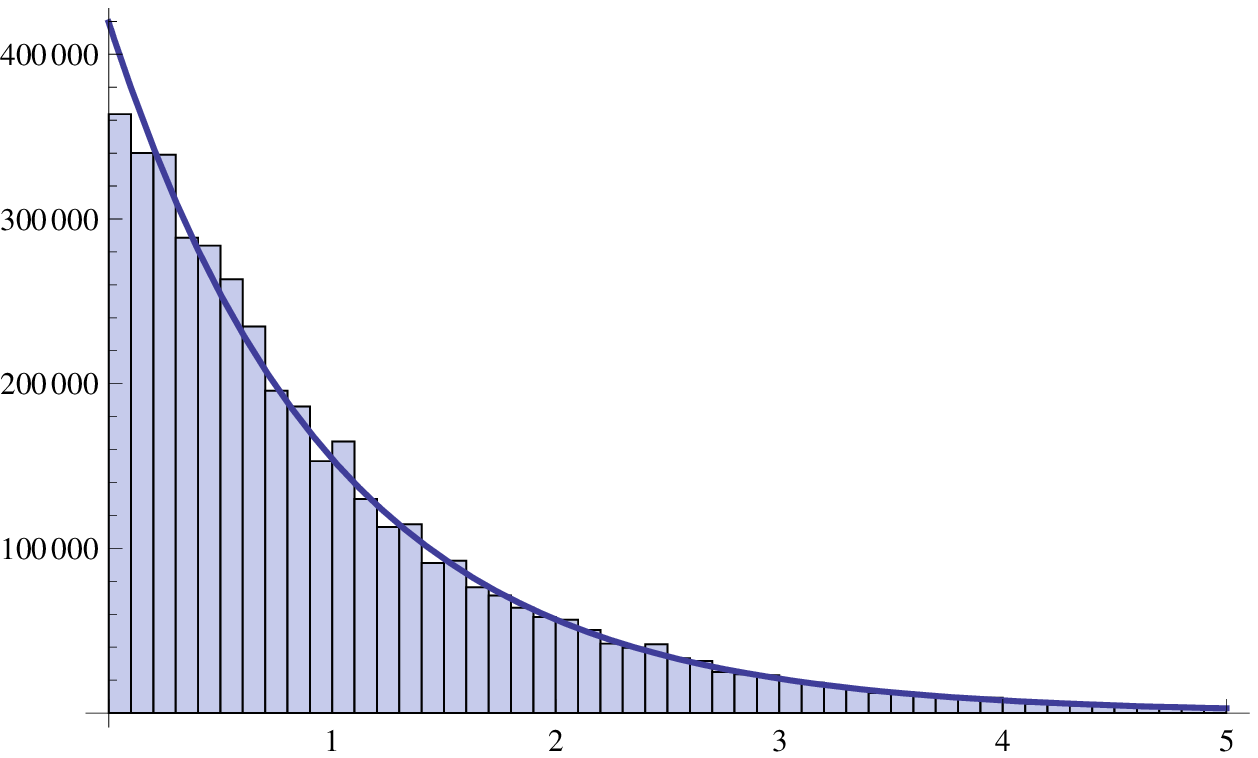} \ \ & \ \
        \epsfxsize=2.5in
      \epsffile{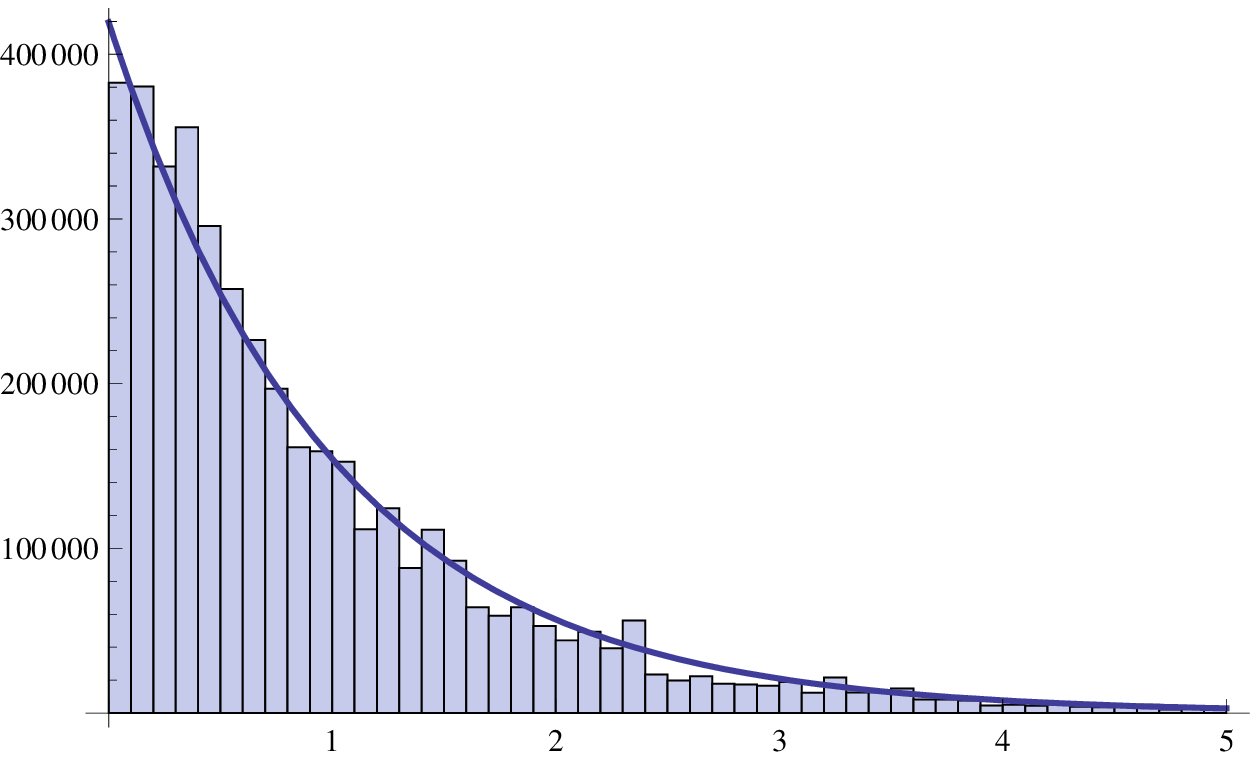} \\[1.cm]
 \mbox{rescaled} \ \ & \ \ \mbox{non-rescaled}
\end{array}$
\end{center}
\caption{Spacings for a ``random'' $\lam=.6429,\ N=22$.}
\end{figure}

Similar results (though, admittedly, some showing less accurate fit) were obtained for more than 20 values of $\lam$
in the range $[0.62,0.78]$. For some values we went up to $N=25$, but this significantly increased the computing time without noticeable advantage.

We also tried $\lam$ close to a Pisot reciprocal and to Garsia reciprocal. In Figure 5 we took $\lam = 0.6518$ which is close to the Pisot reciprocal
$0.651822\ldots$; we still see ``exponential-like'' distribution, but there is definite ``attraction'' visible. For $\lam = 0.652$ this effect  was no longer
felt. In Figure 6 we took $\lam = 0.70710678$ which is within about $10^{-9}$ to $2^{-1/2}$, a Garsia reciprocal.  It shows a very dramatic ``repulsion.''
This effect was not felt for $\lam = 0.7071$.

\begin{figure}[ht]
\begin{center}
$\begin{array}{cc}
\multicolumn{1}{l}{} &
        \multicolumn{1}{l}{} \\
\epsfxsize=3.5in
\epsffile{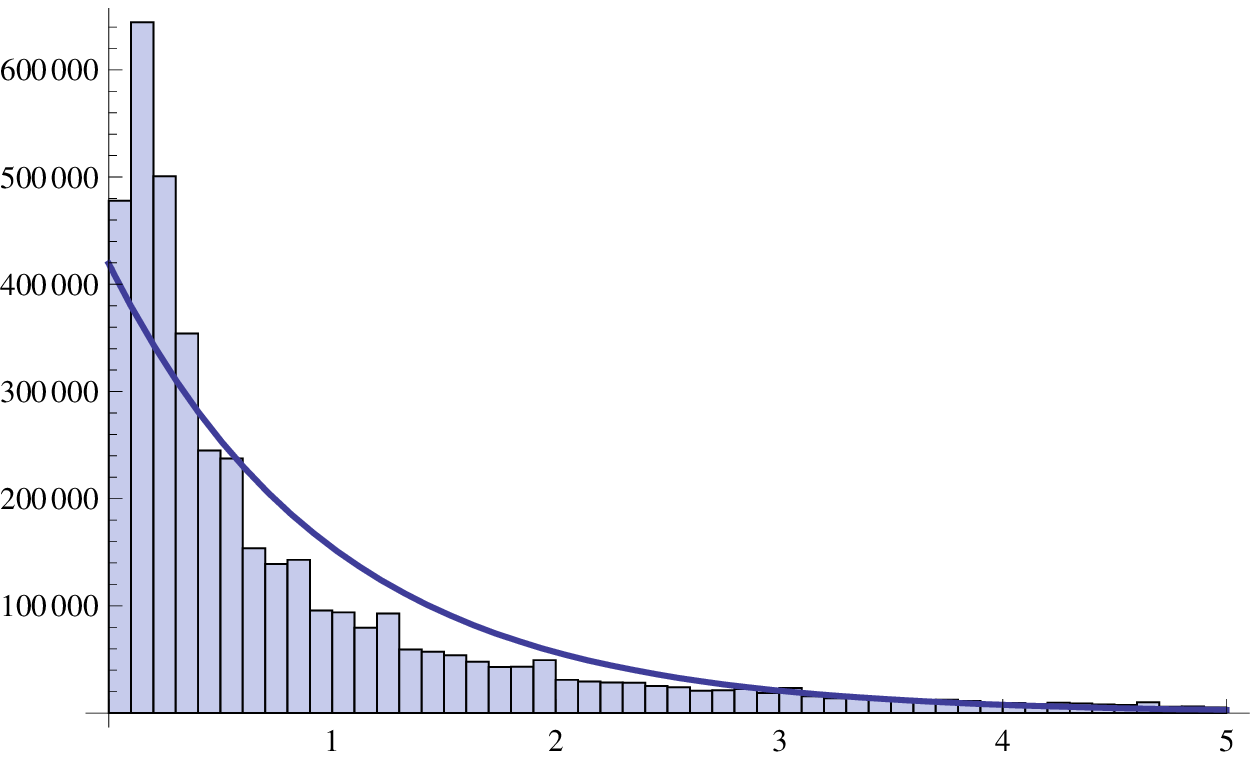} \ \ & \ \
        \epsfxsize=2.5in
      \epsffile{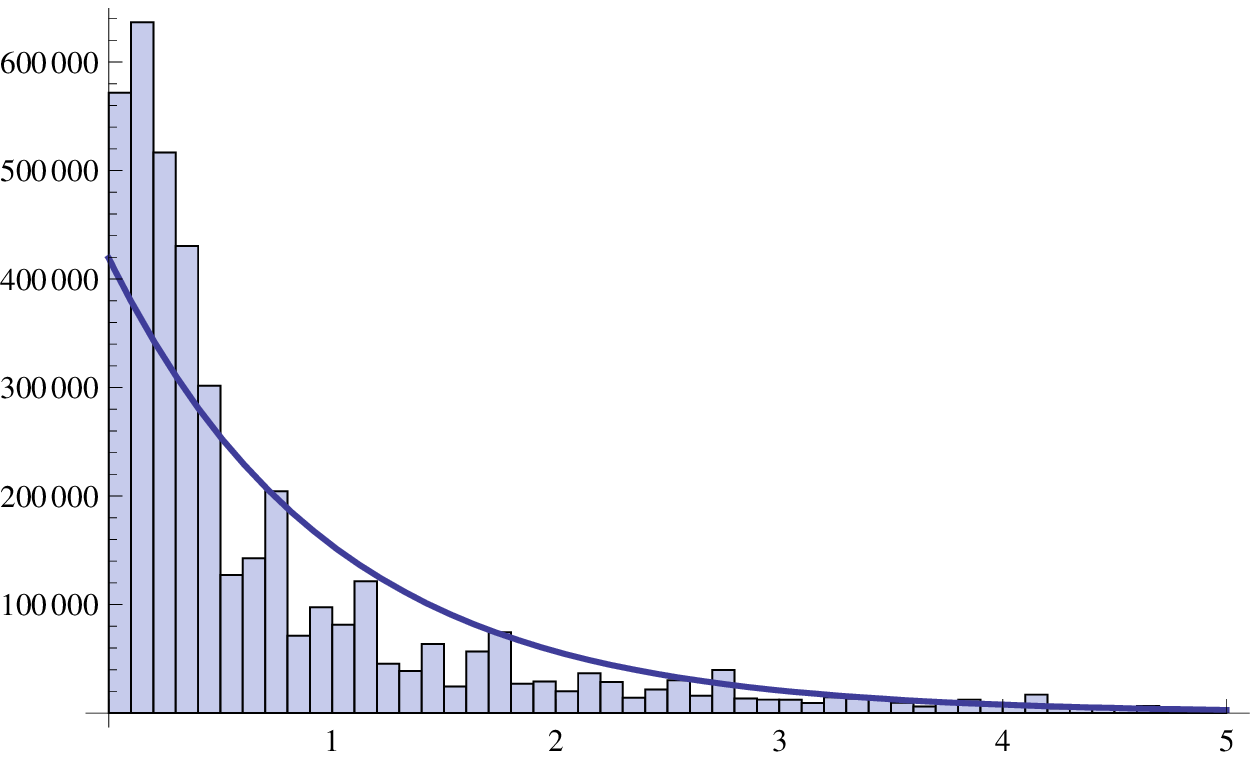} \\[1.cm]
 \mbox{rescaled} \ \ & \ \ \mbox{non-rescaled}
\end{array}$
\end{center}
\caption{Near Pisot: $\lam=.6518,\ N=22$.}
\end{figure}

\begin{figure}[ht]
\begin{center}
$\begin{array}{cc}
\multicolumn{1}{l}{} &
        \multicolumn{1}{l}{} \\
\epsfxsize=3.5in
\epsffile{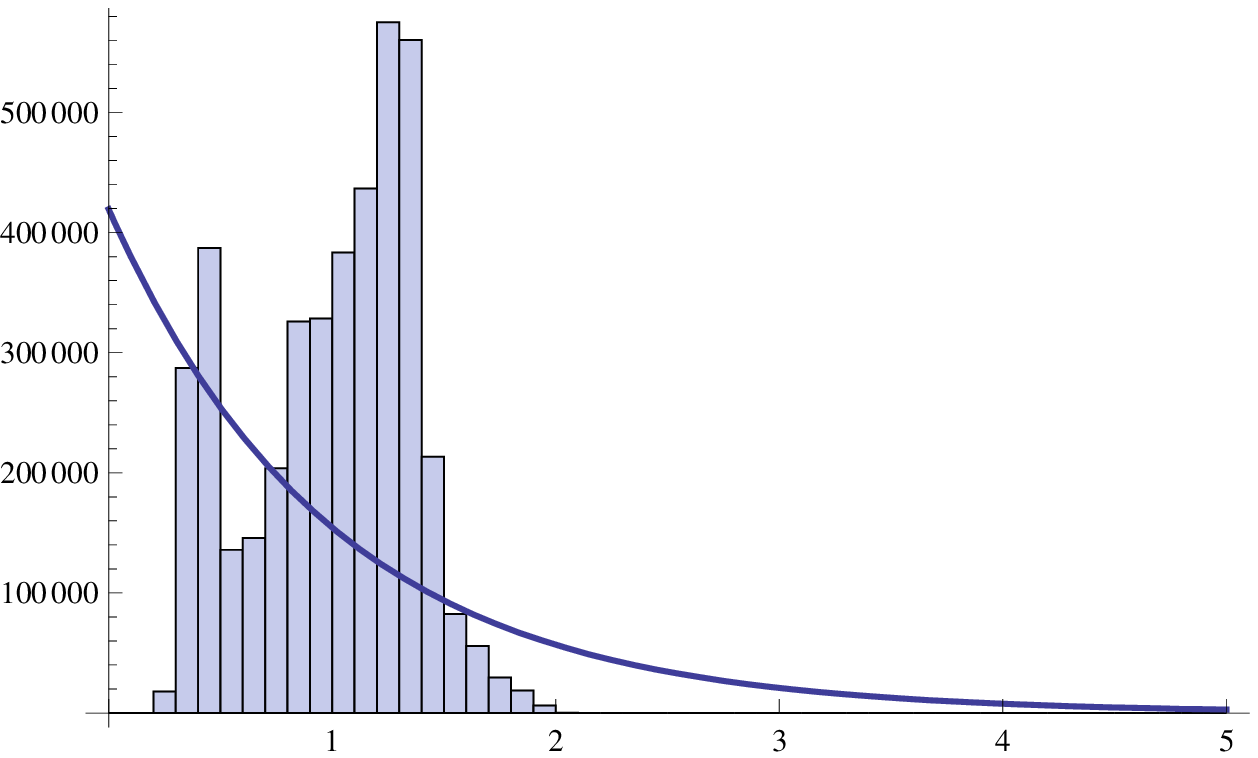} \ \ & \ \
        \epsfxsize=2.5in
      \epsffile{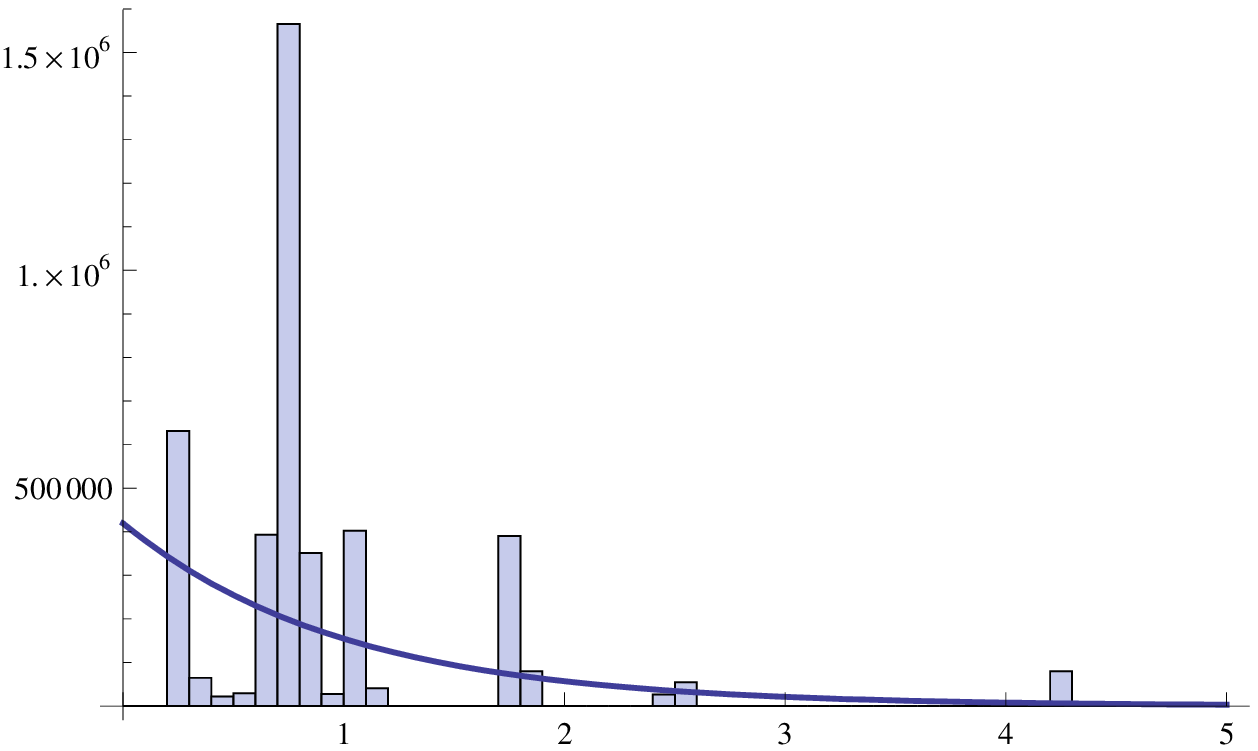} \\[1.cm]
 \mbox{rescaled} \ \ & \ \ \mbox{non-rescaled}
\end{array}$
\end{center}
\caption{Near Garsia: $\lam=.70710678,\ N=22$.}
\end{figure}

To conclude, we provide the {\em Mathematica 6.0} program used to create these figures. We use $\lam=$\verb+lam+$=0.6429$,
$N=22$;
$\ell=$\verb+ell+$=1$ means that we consider nearest-neighbor spacings. We then just change these parameters. 
\verb+Spac+ and \verb+SpacR+ produce the histograms of spacings for non-rescaled and rescaled points, respectively.

\begin{verbatim}
Needs["Histograms`"]
a = 1.5*Sqrt[2.] + 2.; b = Sqrt[2.] - 1.
F[x_] := Which[x < b, a*x^2/2., x > 1 - b,
1 - a (1 - x)^2/2., True, a*b^2/2 + a*b*(x - b)]
lam = 0.6429; N= 22; ell=1; Clear[g]; g[A_] := Union[lam A, lam A + 1]
Tru = Nest[g, {0}, N] (1 - lam);
Vru = Drop[RotateLeft[Pru,ell] - Pru, -ell] 2^N;
Spac = Histogram[Vru, HistogramRange -> {0, 5*ell},
HistogramCategories -> 50, ApproximateIntervals -> False]
Pru = Table[F[Tru[[i]]], {i, 1, 2^N}];
Vuru = Drop[RotateLeft[Pru,ell] - Pru, -1] 2^N;
SpacR = Histogram[Vuru, HistogramRange -> {0, 5*ell},
HistogramCategories -> 50, ApproximateIntervals -> False]
\end{verbatim}

We compared the histograms with appropriately rescaled Poisson distribution  as follows:

\begin{verbatim}
Pois[s_,ell_,N_]:=.1*ell*2^N*s^(ell-1) Exp[-s]/(ell-1)!
Plot[Pois[s,ell,n], {s, 0, 5*ell}, PlotRange -> All]
\end{verbatim}

\medskip

{\bf Acknowledgment.} We are grateful to Rick Kenyon for his ideas at the initial stage of this work. Thanks also to Yury Makarychev for some help with
simulations.

\bibliographystyle{amsplain}

\end{document}